# Instant Evaluation and Demystification of $\zeta(n), L(n,\chi)$ that Euler, Ramanujan Missed-II

VIVEK V. RANE

THE INSTITUTE OF SCIENCE

15, MADAME CAMA ROAD,

MUMBAI-400 032

INDIA

e-mail address : v_v_rane@yahoo.co.in

**Abstract :** We give a formal proof of the power series expression in $\alpha$ of $\zeta_k(s,\alpha) = \zeta(s, k+\alpha)$ For $|\alpha| < k$, where $\zeta(s,\alpha)$ is the Hurwitz zeta function and $k \geq 0$ is an integer. We also show that for an integer $r \geq 0$, the power series expression of $\frac{\partial^r}{\partial s^r}\zeta_1(s,\alpha) = \zeta_1^{(r)}(s,\alpha)$ can be obtained from that of $\zeta_1(s,\alpha)$ by term-by-term differentiation r times . Consequently, from the power series expression of $\zeta'(-m,\alpha)$ (where $m \geq 0$ is an integer ), we obtain rapidly decreasing series for $\zeta(2m+1)$, where $\zeta(s) = \zeta(s,1)$. We obtain analogous series expressions in the case of Dirichlet L-series. Also using the power series ( in $\alpha$ ) of $\zeta_k(s,\alpha)$, we give a unified proof of the well-known classical identities on series involving $\zeta(s)$.

# Instant Evaluation and Demystification of $\zeta(n), L(n,\chi)$ that Euler, Ramanujan Missed-II

VIVEK. V. RANE

THE INSTITUTE OF SCIENCE

15, MADAM CAMA ROAD,

MUMBAI-400 032

INDIA

e-mail address : v_v_rane@yahoo.co.in

Let $s = \sigma + it$ be a complex variable, where $\sigma$ and $t$ are real. For complex $\alpha \neq 0, -1, -2, \ldots$, let $\zeta(s,\alpha)$ be the Hurwitz zeta function defined by $\zeta(s,\alpha) = \sum_{n \geq 0}(n+\alpha)^{-s}$ for $\sigma > 1$ and its analytic continuation. For integral $k \geq 1$, let $\zeta_k(s,\alpha) = \sum_{n \geq k}(n+\alpha)^{-s}$ for $\sigma > 1$, and its analytic continuation. Thus $\zeta_k(s,\alpha) = \zeta(s, k+\alpha) = \zeta(s,\alpha) - \sum_{0 \leq n \leq k-1}(n+\alpha)^{-s}$.

In author [3], we have dealt with the case $k=1$ and have shown that $\zeta_1(s,\alpha)$ is an analytic function of $\alpha$ for a fixed complex number $s$ in the unit disc $|\alpha|<1$ of complex plane. We had also given the power series of $\zeta_1(s,\alpha)$ in the unit disc $|\alpha|<1$ in the form

$$\zeta_1(s,\alpha) - \zeta(s) = \sum_{n \geq 1}\alpha^n \frac{(-1)^n}{n!}s(s+1)\ldots\ldots(s+n-1)\zeta(s+n)$$

By choosing $\alpha = \pm \frac{1}{2}$ and $\alpha = 1$ there, we have obtained three identities involving $\zeta(s)$. However, Ramanujan [2] has already given the power series expression of

$\zeta(s, 1-\alpha)$ without proof and without elaborating on it.

In this paper, our Theorem 1 below extends this result to $\zeta_k(s,\alpha)$ for any integer $k \geq 1$ in the disc $|\alpha| < k$. As a consequence, we shall show that various identities involving $\zeta(s)$, are simple particular cases of the power series expression of $\zeta_k(s,\alpha)$. In fact, we can generate many more interesting identities involving $\zeta(s)$ by choosing various values of integer $k \geq 1$ and for each k by choosing integer values of $\alpha$ with $|\alpha| < k$ or by choosing $\alpha = \pm \frac{1}{2}$. This will be clear from the following examples. There are some apparently unconnected identities involving $\zeta(s)$ which are almost of classical nature, namely

I) $\zeta(s) - \frac{1}{s-1} = 1 - \sum_{n \geq 1} \frac{s(s+1)........(s+n-1)(\zeta(s+n)-1)}{(n+1)!}$ .

II) $(1 - 2^{1-s})\zeta(s) = \sum_{n \geq 1} \frac{s(s+1)........(s+n-1)\zeta(s+n)}{2^{s+n} n!}$ .

III) $\sum_{n \geq 1} \frac{s(s+1)........(s+n-1)}{n!}(\zeta(s+n) - 1) = 1$ .

IV) $\sum_{n \geq 1} \frac{(-1)^{n-1} s(s+1)........(s+n-1)}{n!}(\zeta(s+n) - 1) = 2^{-s}$ .

**Remark :** Letting $s \to 1$ in identity I), we have $\gamma = 1 - \sum_{n \geq 2} \frac{\zeta(n)}{n}$, where $\gamma$ is Euler's constant.

The first two identities have been stated in Titchmarsh's book [ 5 ], of which the first one has been reproduced there from Landau's Handbuch and there, the second one has been attributed to Ramaswami. We shall show that the above identities are all particular cases of our Theorem 1 below. Our Theorem 1 below shows that for a fixed complex number s, $\zeta(s,\alpha) - \zeta(s)$ is an analytic function of the complex variable $\alpha$ in the whole complex plane with singularities at non-positive integers. Our Theorem 2 below shows that for any integer $r \geq 0$, the power series of $\zeta_1^{(r)}(s,\alpha) = \frac{\partial^r}{\partial s^r}\zeta_1(s,\alpha)$ can be obtained by term-by-term differentiation r times, of the power series of $\zeta_1(s,\alpha)$ with respect to s. More generally, this result can be obtained for $\frac{\partial^r}{\partial s^r}\zeta_k(s,\alpha)$, where $k \geq 1$ is an integer. Our Theorem 3 below evaluates $\zeta'(-m,\alpha)$ for any integral $m \geq 0$, in the same fashion as we had obtained $\zeta(-m,\alpha)$ from the power series of $\zeta_1(s,\alpha)$ in author [4]. Here as elsewhere, dash denotes differentiation with respect to s. More generally for integral $r \geq 0, \zeta^{(r)}(s,\alpha)$ will stand for $\frac{\partial^r}{\partial s^r}\zeta(s,\alpha)$. Incidentally, S. Ramanujan has obtained power series expression for $\zeta_1'(-m,\alpha)$. (See entry 28(b) of Ramanujan's Notebooks Part-I by B.C. Berndt [1]. But our approach is completely different and we obtain Theorem 3 as a consequence of Theorem 2, which is important.

Next, we state our Theorems and their corollaries.

**Theorem 1**: We have for an integer $k \geq 1$, and for a fixed complex number s,

$$\zeta_k(s,\alpha) - \zeta_k(s) = \sum_{n \geq 1} \frac{(-1)^n s(s+1)\ldots(s+n-1)}{n!} \zeta_k(s+n)\alpha^n$$ in the disc $|\alpha| < k$ of the complex plane, where

$$\zeta_k(s) = \zeta(s) - \sum_{1 \leq n \leq k-1} n^{-s}.$$

**Corollary**: We have for any integral $m \geq 0$ and for a fixed integral $k \geq 1$ and for a complex $\alpha$ with $|\alpha| < k$, $\zeta(-m,\alpha) = \sum_{0 \leq n \leq k-1}(n+\alpha)^m + \sum_{\ell=0}^{m}\binom{m}{\ell}\left(\zeta(-\ell) - \sum_{1 \leq n \leq k-1} n^\ell\right)\alpha^{m-\ell} - \frac{\alpha^{m+1}}{m+1}$

**Remarks**: 1) Identity I) follows by taking k=2, $\alpha = -1$ and by replacing s with (s-1).

This gives $\zeta_2(s-1,-1) - \zeta_2(s-1) = \sum_{n \geq 1}\frac{(s-1)s\ldots(s+n-2)}{n!}\zeta_2(s-1+n)$

Note that $\zeta_2(s,-1) = \zeta(s)$ and. $\zeta_2(s) = \zeta(s) - 1$.

This gives $1=(s-1)\zeta_2(s)+\sum_{n\geq 2}\frac{(s-1)s...(s+n-2)}{n!}(\zeta(s+n-1)-1)$.

That is, $\frac{1}{s-1}=\zeta_2(s)+\sum_{n\geq 2}\frac{s(s+1)...(s+n-2)}{n!}(\zeta(s+n-1)-1)$.

That is $\frac{1}{s-1}=(\zeta(s)-1)+\sum_{n\geq 1}\frac{s(s+1)...(s+n-1)}{(n+1)!}(\zeta(s+n)-1)$.

2) Similarly identity II) can be obtained by considering $\zeta_1(s,-½)$ and by observing

$$\zeta_1(s,-½)=(2^s-1)\zeta(s).$$

3) Identity III) can be obtained by considering $\zeta_1(s,1)$ and by observing

$$\zeta_1(s,1)=\zeta(s)-1 \quad \text{and} \quad \zeta_1(s)=\zeta(s).$$

4) Identity IV) can be obtained by considering $\zeta_2(s,1)$ and by noting

$$\zeta_2(s,1)=\zeta(s)-1-2^{-s} \quad \text{and} \quad \zeta_2(s)=\zeta(s)-1.$$

Next, we state our Theorem 2.

**Theorem 2 :** For an integer $r\geq 0$ and for a fixed complex s, $\frac{\partial^r}{\partial s^r}\zeta_1(s,\alpha)$ is an analytic function of $\alpha$ in the unit disc $|\alpha|<1$ of the complex plane, where it has a power series expression namely

$$\frac{\partial^r}{\partial s^r}\zeta_1(s,\alpha)-\frac{\partial^r}{\partial s^r}\zeta(s)=\sum_{n\geq 1}\frac{\partial^r}{\partial s^r}\{s(s+1)....(s+n-1)\zeta(s+n)\}\frac{(-\alpha)^n}{n!}.$$

**Note :** For r = 2, we get the following corollary.

**Corollary :** We have $\zeta''(0,\alpha)=\log^2\alpha+\zeta''(0)+2\,\gamma_1\,\alpha+\sum_{n\geq 2}\frac{(-\alpha)^n}{n}\left(2\zeta'(n)+\zeta(n)\sum_{k=1}^{n-1}\frac{1}{k}\right)$,

where $\zeta(s)=\frac{1}{s-1}+\gamma+\sum_{n\geq 1}\gamma_n(s-1)^n$ for $s\neq 1$.

As a consequence of Theorem 2, we get

**Theorem 3** : For $|\alpha|<1$ and for any integer $m \geq 0,$ we have

$$\zeta_1'(-m,\alpha)=\zeta'(-m)+\sum_{k=0}^{m-1}\alpha^{m-k}\binom{m}{k}\left(\zeta'(-k)-\zeta(-k)\sum_{j=k+1}^{m}\tfrac{1}{j}\right)+\tfrac{\alpha^{m+1}}{m+1}\left(\sum_{j=1}^{m}\tfrac{1}{j}-\gamma\right)+\sum_{k\geq 2}\tfrac{(-1)^k(k-1)!m!}{(m+k)!}\zeta(k)\alpha^{m+k},$$

where empty sum stands for zero and $\gamma$ stands for Euler's constant.

**Corollary:** We have $\log\Gamma(\alpha+1)=-\gamma\alpha+\sum_{n\geq 2}(-1)^n\tfrac{\zeta(n)\alpha^n}{n}$ for any complex number α with $|\alpha|<1$.

Next, we give the proof of our Theorems.

**Proof of Theorem 1** : Our proof follows closely the proof of Theorem 1 of author [3], where we have dealt with $\zeta_k(s,\alpha)$ for k=1. For $\sigma>1$ and for $|\alpha|<k,$ consider

$$\zeta_k(s,\alpha)=\sum_{n\geq k}(n+\alpha)^{-s}=\sum_{n\geq k}n^{-s}(1+\tfrac{\alpha}{n})^{-s}=\sum_{n\geq k}n^{-s}\left\{1-s\tfrac{\alpha}{n}+\tfrac{s(s+1)}{2!}\tfrac{\alpha^2}{n^2}-\ldots\ldots\right\}$$

$$=\sum_{n\geq k}n^{-s}-\alpha s\sum_{n\geq k}n^{-s-1}+\tfrac{\alpha^2 s(s+1)}{2!}\sum_{n\geq k}n^{-s-2}-\ldots\ldots$$

$$=\zeta_k(s)-\alpha s\zeta_k(s+1)+\tfrac{s(s+1)\alpha^2}{2!}\zeta_k(s+2)-\ldots\ldots=\sum_{n\geq 0}b_{n,k}(s)\alpha^n, \text{ say}.$$

Next we show the series $\sum_{n\geq o}b_{n,k}(s)\alpha^n$ converges absolutely and uniformly in every compact subset of the disc $|\alpha|<k$ when $\sigma>1$. Hence consider $\beta$ such that $o<|\alpha|\leq \beta<k$.

For $\sigma>1$, $|\sum_{n\geq o}b_{n,k}(s)\alpha^n|\leq\sum_{n\geq o}|b_{n,k}(s)||\alpha|^n\leq\sum_{n\geq o}|b_{n,k}(s)|\beta^n$

Let $\beta'=\tfrac{\beta}{k}$ so that $\beta'<1$. Thus, we have $|\sum_{n\geq 0}b_{n,k}(s)\alpha^n|\leq\sum_{n\geq 0}k^n|b_{n,k}(s)|\beta'^n.$

Consider

$$k^n|b_{n,k}(s)|\leq\tfrac{k^n s(s+1)\ldots(s+n-1)|\zeta_k(s+n)|}{n!}\leq\tfrac{|s|(|s|+1)\ldots(|s|+n-1)k^n}{n!}\sum_{m\geq k}m^{-\sigma-n}.$$

$$\leq\tfrac{|s|(|s|+)\ldots(|s|+n-1)}{n!}\zeta(\sigma)$$

Thus $|\sum_{n\geq 0} b_{n,k}(s)\alpha^n | \leq \zeta(\sigma) \sum_{n\geq 0} \frac{|s|(|s|+1)\dots(|s|+n-1)\beta'^n}{n!} \leq \zeta(\sigma)(1-\beta')^{-|s|}$.

Thus $\sum_{n\geq 0} b_{n,k}(s)\alpha^n$ defines an analytic function of $\alpha$ for a fixed s with $\sigma > 1$, where

$\zeta_k(s,\alpha) = \sum_{n\geq 0} b_{n,k}(s)\alpha^n$. Consequently, we have $\zeta_k(s,\alpha) = \sum_{n\geq 0} b_{n,k}(s)\alpha^n$ for any complex number s.

**Proof of Theorem 2 :** We have for $\sigma > 1$, $\zeta_1^{(r)}(s,\alpha) = (-1)^r \sum_{n\geq 1}(n+\alpha)^{-s} \log^r(n+\alpha)$

Thus for $\sigma > 1, (-1)^r \zeta_1^r(s,\alpha) = \sum_{n\geq 1}(n+\alpha)^{-s} \log^r(n+\alpha) = \sum_{n\geq 1} n^{-s}(1+\frac{\alpha}{n})^{-s}(\log n(1+\frac{\alpha}{n}))^r$,

$= \sum_{n\geq 1} n^{-s}(1+\frac{\alpha}{n})^{-s} (\log n + \log(1+\frac{\alpha}{n}))^r = \sum_{n\geq 1} n^{-s}(1+\frac{\alpha}{n})^{-s} \cdot \sum_{\ell=0}^{r}\binom{r}{\ell} \log^{r-\ell} n \cdot \log^\ell(1+\frac{\alpha}{n})$

$= \sum_{\ell=0}^{r}\binom{r}{\ell} \cdot \sum_{n\geq 1} n^{-s}(1+\frac{\alpha}{n})^{-s} \log^{r-\ell} n \cdot \log^\ell(1+\frac{\alpha}{n})$.

Consider $\sum_{n\geq 1} n^{-s} \log^{r-\ell} n \cdot (1+\frac{\alpha}{n})^{-s} \cdot (\log(1+\frac{\alpha}{n}))^\ell = \sum_{n\geq 1} n^{-s} \log^{r-\ell} n \cdot \{1 - s\frac{\alpha}{n} + \frac{s(s+1)}{2!}(\frac{\alpha}{n})^2 - \dots\}$

$\{\frac{\alpha}{n} - \frac{1}{2}(\frac{\alpha}{n})^2 + \frac{1}{3}(\frac{\alpha}{n})^3 - \dots\}^\ell = \sum_{n\geq 1} n^{-s} \log^{r-\ell} n \cdot \left(\sum_{i\geq 0} a_i(s)(\frac{\alpha}{n})^i\right)\left(\sum_{j\geq 1} b_j(\frac{\alpha}{n})^j\right)^\ell$ say.

$= \sum_{n\geq 1} n^{-s} \log^{r-\ell} n \cdot \left(\sum_{i\geq 0} a_i(s)(\frac{\alpha}{n})^i\right)\left(\sum_{j\geq \ell} c_j(\frac{\alpha}{n})^j\right)$, say. $= \sum_{n\geq 1} n^{-s} \log^{r-\ell} n \cdot \sum_{k\geq r} d_k(s) \cdot (\frac{\alpha}{n})^k$, say.

$= \sum_{k\geq r} d_k(s)\alpha^k \left(\sum_{n\geq 1} n^{-s-k} \log^{r-\ell} n\right) = (-1)^{r-\ell} \sum_{k\geq \ell} d_k(s) \zeta^{(r-\ell)}(s+k)\alpha^k$.

Thus $\zeta_1^{(r)}(s,\alpha) = \sum_{\ell=0}^{r}(-1)^\ell \binom{r}{\ell} \sum_{k\geq \ell} d_k(s)\zeta^{(r-\ell)}(s+k)\alpha^k$.

Next, we show that for a fixed complex number s, the series for $\zeta_1^{(r)}(s,\alpha)$ converges uniformly on every compact subset of the unit disc $|\alpha|<1$. For this, we show that for each $\ell$ with $0\leq \ell \leq r$, $\sum_{k\geq r} d_k(s)\zeta^{(r-\ell)}(s+k)\alpha^k$ converges uniformly on every compact subset of the disc $|\alpha|<1$. Consider the compact subset on which $|\alpha|\leq \beta <1$. For simplicity, we consider a fixed s with Re s > 1. Consider

$$\left|\sum_{k\geq \ell} d_k(s)\zeta^{(r-\ell)}(s+k)\alpha^k\right| \leq \sum_{k\geq \ell}|d_k(s)| \ |\zeta^{(r-\ell)}(\sigma)\alpha^k| \leq |\zeta^{(r-\ell)}(\sigma)|\sum_{k\geq \ell}|d_k(s)|\alpha^k.$$

Next consider $\sum_{k\geq \ell}|d_k(s)|\alpha^k = \sum_{k\geq \ell}\left|\sum_{i+j=k}a_i(s)\cdot c_j\right|\alpha^k \leq \sum_{k\geq \ell}\left(\sum_{i+j=k}|a_i(s)|\cdot |c_j|\right)\alpha^k$

$$\leq \left(\sum_i |a_i(s)|\alpha^i\right)\left(\sum_j |c_j|\alpha^j\right) \leq \left(\sum_i |a_i(s)|\alpha^i\right)\left(\sum_{j_1+j_2+..+j_\ell=j}|b_{j_1}b_{j_2}...b_{j_\ell}|\right)\alpha^j$$

$$\leq \left(\sum_i|a_i(s)|\alpha^i\right)\left(\sum_{j_1}|b_{j_1}|\alpha^{j_1}\right)\left(\sum_{j_2}|b_{j_2}|\alpha^{j_2}\right)........\left(\sum_{j_\ell}|b_{j_\ell}|\alpha^{j_\ell}\right)\leq \left(\sum_i|a_i(s)|\alpha^i\right)\left(\sum_j|b_j|\alpha^j\right)^\ell$$

Note that $|a_i(s)|=\left|\frac{(-1)^i s(s+1).......(s+i-1)}{i!}\right|\leq \frac{|s|(|s|+1).......(|s|+i-1)}{i!}$ so that

$$\sum_i |a_i(s)|\alpha^i \leq \sum_i \frac{|s|(|s|+1).......(|s|+i-1)}{i!}\beta^i \leq (1-\beta)^{-|s|}.$$

Similarly $|b_j|\leq \frac{1}{j}$ so that $\sum_{j\geq 1}|b_j|\alpha^j \leq \sum_{j\geq 1}\frac{\beta^j}{j}\leq -\log(1-\beta)$.

Thus $\left|\sum_{k\geq \ell}d_k(s)\zeta^{(r-\ell)}(s+k)\alpha^k\right|\leq (1-\beta)^{-|s|}(-\log(1-\beta))^r |\zeta^{(r-\ell)}(\sigma)|$

Thus for a fixed s, $\zeta_1^{(r)}(s,\alpha)$ is an analytic function of $\alpha$ in the unit disc $|\alpha|<1$ and hence

$$\zeta_1^{(r)}(s,\alpha) = \sum_{k \geq 0} \frac{A_k(s)\alpha^k}{k!}, \text{ where } A_k(s) = \left(\frac{\partial^k}{\partial \alpha^k} \frac{\partial^r}{\partial s^r} \zeta_1(s,\alpha)\right)_{\alpha=0}$$

We shall show for any complex number s, that $\left(\frac{\partial^k}{\partial \alpha^k} \frac{\partial^r}{\partial s^r} \zeta_1(s,\alpha)\right)_{\alpha=0} = \frac{\partial^r}{\partial s^r} \left(\frac{\partial^k}{\partial \alpha^k} \zeta_1(s,\alpha)\right)_{\alpha=0}$

This will prove our contention that the power series for $\zeta_1^{(r)}(s,\alpha)$ can be obtained by r times term-by-term differentiation (with respect to s) of the power series

$$\zeta_1(s,\alpha) = \zeta_1(s) + \sum_{n \geq 1} \frac{(-1)}{k!} s(s+1)\ldots(s+n-1)\zeta(s+n)\alpha^n$$

Since $\left(\frac{\partial^k}{\partial \alpha^k} \frac{\partial^r}{\partial s^r} (n+\alpha)^{-s}\right)_{\alpha=0} = \frac{\partial^r}{\partial s^r}\left(\frac{\partial^k}{\partial s^k}(n+\alpha)^{-s}\right)_{\alpha=0}$ for any integer $n \geq 1$,

we have for $\sigma > 1$, $\left(\frac{\partial^k}{\partial \alpha^k} \frac{\partial^r}{\partial s^r} \zeta_1(s,\alpha)\right)_{\alpha=0} = \frac{\partial^r}{\partial s^r}\left(\frac{\partial^k}{\partial \alpha^k} \zeta_1(s,\alpha)\right)_{\alpha=0}$

Note that either side is an analytic function of s and thus being equal for $\sigma > 1$, the equality follows for any complex $s \neq 1$..

**Proof of Theorem 3 :** $\zeta_1(s,\alpha) = \zeta(s) + \sum_{n \geq 1} \frac{s(s+1)\ldots(s+n-1)\zeta(s+n)}{n!}(-\alpha)^n$ Hence

$$\zeta_1'(s,\alpha) = \zeta'(s) + \sum_{n \geq 1} \left(s(s+1)\ldots(s+n-1)\zeta(s+n)\right)' \frac{(-\alpha)^n}{n!}$$

where dash over the bracket as elsewhere denotes differentiation with respect to s. Next we obtain

$\zeta_1'(-m,\alpha)$, where m≥1 is an integer.

Now $\zeta_1'(-m,\alpha) = \zeta'(-m) + \left(\sum_{1 \leq n \leq m} + \sum_{n=m+1} + \sum_{n \geq m+2}\right) \left(s(s+1)\ldots(s+n-1)\zeta(s+n)\right)'_{s=-m} \frac{(-\alpha)^n}{n!}$

Consider $\left(s(s+1)\ldots(s+n-1)\zeta(s+n)\right)'_{s=-m}$ for n=m+1.

This is equal to $\left(s(s+1)\ldots(s+m)\zeta(s+m+1)\right)'_{s=-m} = \left((s(s+1)\ldots(s+m-1))'(s+m)\zeta(s+m+1)\right)_{s=-m}$

$+ \left(s(s+1)\ldots(s+m-1)\right)_{s=-m} \cdot \left((s+m)\zeta(s+m+1)\right)'_{s=-m}$

Note for $s = -m$, $(s+m)\zeta(s+m+1) = \lim_{w \to o} w\zeta(w+1) = 1$ and
$((s+m)\zeta(s+m+1))' = \lim_{w \to o}\left(\frac{d}{dw}(w\zeta(w+1))\right) = \gamma$, the Euler's constant.

Thus
$$(s(s+1)....(s+m)\zeta(s+m+1))'_{s=-m} = (s(s+1)...(s+m-1))'+\gamma(s(s+1)...(s+m-1))_{s=-m}$$

$$= \left(\prod_{i=0}^{m-1}(s+i)\sum_{j=0}^{m-1}\frac{1}{s+j}\right)_{s=-m} + \gamma(-1)^m m! = (-1)^m m!\left(\gamma - \sum_{j=1}^{m}\frac{1}{j}\right)$$

Thus $(s(s+1)...(s+m)\zeta(s+m+1))'_{s=-m} \cdot \frac{(-\alpha)^{m+1}}{(m+1)!} = \frac{\alpha^{m+1}}{m+1}\left(\sum_{j=1}^{m}\frac{1}{j} - \gamma\right)$

Next consider $\sum_{n \geq m+2}(s(s+1)...(s+n-1)\zeta(s+n))'_{s=-m}\frac{(-\alpha)^n}{n!}$

Now for $n \geq m+2$, $(s(s+1)...(s+n-1)\zeta(s+n))'_{s=-m} = ((s(s+1)...(s+n-1)\zeta'(s+n))_{s=-m}$

$+((s(s+1)...(s+n-1))' \cdot \zeta(s+n))_{n=-m} = ((s(s+1)...(s+n-1))'\zeta(s+n))_{n=-m}$,

since $(s(s+1)...(s+n-1)\zeta'(s+n))_{s=-m}$ vanishes as it has a vanishing factor s+m for each n≥ m+2.

Thus for each n≥ m+2,
$$\sum_{n \geq m+2}\frac{(-\alpha)^n}{n!}(s(s+1)....(s+n-1)\zeta(s+n))_{n=-m}$$

$$= \sum_{n \geq m+2}\left\{(s(s+1)...(s+n-1))'\zeta(s+n)\right\}_{s=-m} \cdot \frac{(-\alpha)^n}{n!}$$

$$= \sum_{n \geq m+2}\frac{(-\alpha)^n}{n!}\left(\prod_{i=0}^{m-1}(i-m)\right)\prod_{i=m+1}^{n-1}(i-m)\zeta(n-m)$$

$$= \sum_{n \geq m+2}(-1)^m m!(n-m-1)!\zeta(n-m)\frac{(-\alpha)^n}{n!} = \sum_{n \geq m+2}(-\alpha)^n\frac{(-1)^m\zeta(n-m)}{(n-m)\binom{n}{m}} = \sum_{n \geq m+2}\frac{(-\alpha)^n(-1)^m\zeta(n-m)}{(n-m)\binom{n}{n-m}} = \sum_{k \geq 2}\frac{(-\alpha)^{m+k}(-1)^m\zeta(k)}{k\binom{m+k}{k}},$$

where we have written n-m=k.

This, in turn $= \sum_{k \geq 2}\frac{(-1)^k(k-1)!m!\zeta(k)\alpha^{m+k}}{(m+k)!}$

: 10 :

Similarly $\sum_{1\leq n\leq m} \frac{(-\alpha)^n}{n!}\left(s(s+1)...(s+n-1)\zeta(s+n)\right)'_{s=-m}$

$= \sum_{1\leq n\leq m} \frac{(-\alpha)^n}{n!}\left\{\prod_{i=0}^{n-1}(s+i)\left(\sum_{j=0}^{n-1}\frac{1}{s+j}\right)\zeta(s+n)+\prod_{i=0}^{n-1}(s+i)\zeta'(s+n)\right\}_{s=-m}$

$= \sum_{1\leq n\leq m} \frac{\alpha^n}{n!}\left\{-\prod_{i=0}^{n-1}(m-i)\sum_{j=0}^{n-1}\frac{1}{m-j}\zeta(n-m)+\prod_{i=0}^{n-1}(m-i)\zeta'(n-m)\right\}$

$= \sum_{1\leq n\leq m} \frac{\alpha^n}{n!}\prod_{i=0}^{n-1}(m-i)\left(\zeta'(n-m)-\zeta(n-m)\sum_{j=0}^{n-1}\frac{1}{m-j}\right)$

$= \sum_{1\leq n\leq m} \alpha^n \binom{m}{n}\left(\zeta'(n-m)-\zeta(n-m)\sum_{j=0}^{n-1}\frac{1}{m-j}\right) = \sum_{1\leq n\leq m} \alpha^n\left\{\binom{m}{m-n}\left(\zeta'(n-m)-\zeta(n-m)\sum_{j=0}^{n-1}\frac{1}{m-j}\right)\right\}$

$= \sum_{k=0}^{m-1} \alpha^{m-k}\left\{\binom{m}{k}\left(\zeta'(-k)-\zeta(-k)\sum_{j=0}^{m-k-1}\frac{1}{m-j}\right)\right\}$, where we have written m-n=k.

This, in turn

$= \sum_{k=0}^{m-1} \alpha^{m-k}\binom{m}{k}\left(\zeta'(-k)-\zeta(-k)\sum_{j=k+1}^{m}\frac{1}{j}\right)$

Since $\zeta(s,\alpha) = \zeta_1(s,\alpha) + \alpha^{-s}$, we get $\zeta'(s,\alpha) = -\alpha^{-s}\log\alpha + \zeta_1'(s,\alpha)$

and the proof follows for m≥1.

For m=o, it is easy to see that

$\zeta'(o,\alpha) = \left\{-\alpha^{-s}\log\alpha + \zeta'(s) - \alpha(\frac{d}{ds}s\zeta(s+1)) + \sum_{n\geq 1}(-\alpha)^n(s+1)(s+2)...(s+n-1)\zeta(s+n)\right\}_{s=0}$

Thus, $\zeta'(o,\alpha) = -\log\alpha + \zeta'(o) - \gamma\alpha + \sum_{n\geq 2}\frac{(-\alpha)^n}{n}\zeta(n)$.

From the fact that $\log\frac{\Gamma(\alpha)}{\sqrt{2\pi}} = \zeta'(o,\alpha)$, where $\Gamma(\alpha)$ denotes gamma function, and the fact that

: 11 :

$\Gamma(\alpha+1) = \alpha\Gamma(\alpha)$, we get the power series for $\log\Gamma(\alpha+1)$ as in the corollary above.

Following our Theorem 3 and Lemma 3 of author [4], we get the following Theorem 4.

**Theorem 4**: I) We have for $m \geq 1$,

$$\zeta(2m+1) = \frac{(-1)^m \pi^{2m} \log 2}{(2^{-2m}-2)(2m)!} + \frac{(-1)^m \cdot 2\pi^{2m}}{(2m)!(2^{-2m}-2)} \left\{ \frac{1}{2}\sum_{j=1}^{2m}\frac{1}{j}\sum_{k=0}^{m-1} 2^{2k}\binom{2m}{2k}\zeta'(-2k) + \sum_{k\geq 1}\frac{(2k-1)!(2m)!}{(2m+2k)!}\zeta(2k)\cdot 2^{-2k} \right\}$$

II) For an integer $q \geq 2$, let $\chi(\bmod q)$ be a primitive, even character and let $m \geq 0$ be an integer. Then we have $L(2m+1, \chi) = \frac{(-1)^m 2^{2m} \pi^{2m}}{(2m)!\tau(\overline{\chi})} \left\{ \log q \sum_a \overline{\chi}(a)a^{2m} - \sum_a \overline{\chi}(a)a^{2m} \log a \right\}$

$$+ \frac{(-1)^m \cdot 2^{2m+1} \pi^{2m}}{(2m)!\tau(\overline{\chi})} \left\{ \sum_{k=0}^{m-1} q^{2k-2m}\binom{2m}{2k}\zeta'(-2k)\sum_a \overline{\chi}(a)a^{2m-k} + \frac{1}{2}\left(\sum_{j=1}^{2m}\frac{1}{j}\right)\sum_a \overline{\chi}(a)a^{2m} \right.$$

$$\left. + \sum_{k=1}^{\infty}\frac{(2k-1)!(2m)!}{(2m+2k)!}\zeta(2k)\cdot q^{-2m-2k}\sum_a \overline{\chi}(a)a^{2m+2k} \right\}$$

III) Let $\chi(\bmod q)$ be an odd, primitive character, where $q \geq 2$ is an integer and let $m \geq 1$ be an integer. Then we have

$$L(2m, \chi) = \frac{i(-1)^{m+1} 2^{2m-1} \pi^{2m-1}}{\tau(\overline{\chi})\cdot(2m-1)!} \cdot \left\{ -\sum_a \overline{\chi}(a)a^{2m-1}\log a + \log q \sum_a \overline{\chi}(a)a^{2m-1} \right\} + \frac{i(-1)^{m+1} 2^{2m} \pi^{2m-1}}{\tau(\overline{\chi})\cdot(2m-1)!} \cdot$$

$$\cdot \left\{ q^{2k+1-2m}\sum_{k=0}^{m-1}\left(\sum_a \overline{\chi}(a)a^{2m-1-2k}\right)\binom{2m-1}{2k}\zeta'(-2k) - \frac{q^{1-2m}}{2}\left(\sum_a \overline{\chi}(a)a^{2m-1}\right)\left(\sum_{j=1}^{2m-1}\frac{1}{j}\right) \right.$$

$$\left. + q^{1-2m-2k}\sum_{k\geq 1}\left(\sum_a \overline{\chi}(a)a^{2m+2k-1}\right)\frac{(2k-1)!\cdot(2m-1)!}{(2m+2k-1)!}\cdot \zeta(2k) \right]$$